\documentclass[11pt]{article}

\usepackage{amscd,amsmath, amssymb, fancyhdr, epsfig, color, tocloft, mathtools}
\usepackage{tikz-cd}

\usepackage[backref=page]{hyperref}
\renewcommand*{\backrefalt}[4]{%
	\ifcase #1 (Not cited.)%
	\or        (Cited on page~#2.)%
	\else      (Cited on pages~#2.)%
	\fi}

\hypersetup{
	colorlinks   = true,
	citecolor    = magenta,
	linkcolor    = blue,
	urlcolor     = magenta	
}

\newcommand{\version}{version 1, \ \ June, 25th, 2026}
\makeatletter
\def\x@arrow{\DOTSB\Relbar}
\def\xlongequalsignfill@{\arrowfill@\x@arrow\Relbar\x@arrow}
\providecommand{\xlongequal}[2][]{%
	\ext@arrow 0099\xlongequalsignfill@{#1}{#2}}
\def\xlongrightarrowfill@{\arrowfill@\relbar\relbar\longrightarrow}
\makeatother

\numberwithin{equation}{section}

\def\eqref#1{(\ref{#1})}

\newcommand{\C}{{\mathbb C}}

\def\1{\sqrt{-1}\:}
\newcommand{\restrict}[1]{{\left|_{{\phantom{|}\!\!}_{#1}}\right.}}
\newcommand{\cntrct}                
{\hspace{2pt}\raisebox{1pt}{\text{$\lrcorner$}}\hspace{2pt}}
\newcommand{\arrow}{{\:\longrightarrow\:}}

\renewcommand{\bar}{\overline}
\renewcommand{\phi}{\varphi}
\renewcommand{\epsilon}{\varepsilon}
\renewcommand{\geq}{\geqslant}

\newcommand{\Reg}{\operatorname{Reg}}

\newcommand{\Def}{\operatorname{Def}}

\renewcommand{\to}{\longrightarrow}

\newcounter{Mycounter}[section]
\newcounter{lemma}[section]
\renewcommand{\thelemma}{{Lemma \thesection.\arabic{lemma}}}
\newcommand{\lemma}{%
	\setcounter{lemma}{\value{Mycounter}}
	\refstepcounter{lemma}
	\stepcounter{Mycounter}
	{\noindent \bf \thelemma:\ }}

\newcounter{claim}[section]
\newcounter{sublemma}[section]
\newcounter{corollary}[section]
\newcounter{theorem}[section]
\renewcommand{\thetheorem}{{Theorem \thesection.\arabic{theorem}}}
\newcommand{\theorem}{%
	\setcounter{theorem}{\value{Mycounter}}
	\refstepcounter{theorem}
	\stepcounter{Mycounter}
	{\noindent \bf \thetheorem:\ }}
	
\newcounter{maintheorem}[section]
\makeatletter 
\newcommand{\maintheorem}{%
	\setcounter{maintheorem}{\value{Mycounter}}%
	\refstepcounter{maintheorem}
	\stepcounter{Mycounter}%
	{\noindent \textbf{Main Theorem:\ }}%
}
\makeatother

\newcounter{conjecture}[section]
\newcounter{proposition}[section]
\renewcommand{\theproposition} {{Proposition \thesection.\arabic{proposition}}}
\newcommand{\proposition}{%
	\setcounter{proposition}{\value{Mycounter}}
	\refstepcounter{proposition}
	\stepcounter{Mycounter}
	{\noindent \bf \theproposition:\ }}

\newcounter{definition}[section]
\renewcommand{\thedefinition} {{Definition~\thesection.\arabic{definition}}}
\newcommand{\definition}{%
	\setcounter{definition}{\value{Mycounter}}
	\refstepcounter{definition}
	\stepcounter{Mycounter}
	{\noindent \bf \thedefinition:\ }}
	
\newcounter{condefinition}[section]
\renewcommand{\thecondefinition} {{Construction and Definition~\thesection.\arabic{condefinition}}}
\newcommand{\condefinition}{%
	\setcounter{condefinition}{\value{Mycounter}}
	\refstepcounter{condefinition}
	\stepcounter{Mycounter}
	{\noindent \bf \thecondefinition:\ }}

\newcounter{example}[section]
\renewcommand{\theexample}{{Example \thesection.\arabic{example}}}
\newcommand{\example}{%
	\setcounter{example}{\value{Mycounter}}
	\refstepcounter{example}
	\stepcounter{Mycounter}
	{\noindent \bf \theexample:\ }}

\newcounter{remark}[section]
\renewcommand{\theremark}{{Remark \thesection.\arabic{remark}}}
\newcommand{\remark}{%
	\setcounter{remark}{\value{Mycounter}}
	\refstepcounter{remark}
	\stepcounter{Mycounter}
	{\noindent \bf \theremark:\ }}

\newcounter{assumption}[section]
\newcounter{problem}[section]
\newcounter{question}[section]
\makeatletter

\@addtoreset{equation}{section}

\makeatother

\def\blacksquare{\hbox{\vrule width 5pt height 5pt depth 0pt}}
\def\endproof{\blacksquare}

\newcommand{\proof}{{\bf Proof: \ }}

\begin{document}
	
\begin{center}
{\Large\bf  Locally Conformally K\"ahler Manifolds\\[.1in] of Algebraic Codimension One}\\[5mm]
{\large
Liviu Ornea\footnote{Partially supported the PNRR-III-C9-2023-I8 grant CF 149/31.07.2023 Conformal Aspects of Geometry and Dynamics.}, 
Misha Verbitsky\footnote{Partially supported by
	FAPERJ SEI-260003/000410/2023 and CNPq - Process 310952/2021-2.}, 
Victor Vuletescu\footnote{Partially supported the PNRR-III-C9-2023-I8 grant CF 149/31.07.2023 Conformal Aspects of Geometry and Dynamics.
\\[1mm]

}\\[4mm]

}
\end{center}

\hfill
	
{\small
\hspace{0.15\linewidth}
\begin{minipage}[t]{0.7\linewidth}
{\bf Abstract} \\
A locally conformally K\"ahler (LCK) manifold
is a manifold $M$ which admits a K\"ahler structure on
its universal cover $\tilde M$, in such a way that the monodromy
acts conformally on $\tilde M$. Let $M$ be an $n$-dimensional
compact LCK manifold of algebraic dimension $n-1$. We prove
that $M$ is bimeromorphic to the total space of an isotrivial
elliptic fibration. Morever, there exists an alteration
of $M$ which dominates bimeromorphically a manifold admitting a free action of an elliptic curve.
\end{minipage} 
}

{\small
	\tableofcontents
}
	
\section{Introduction}

The classification of compact complex manifolds is still a widely open problem. In particular, the classification of compact non-K\"ahler manifolds is a very active field of research. 

The first and most relevent results up to-day are in complex dimension $2$, basically due to Kodaira in the 60's. Nonetheless, a complete classification is still missing, as the understanding of the surfaces in the so-called  {\em Kodaira class} ${\rm VII_0}$ with $b_2>1$,  is still unclear.

However, for non-K\"ahler compact surfaces with algebraic codimension one, the picture is clear (\cite[Proposition 5.1]{_Barth_Peters_Van_de_Ven_}, \cite{_ovv:surf_}): they are all {\em elliptic fibrations}, meaning that for any such surface $S$ there exists a smooth projective curve $B$ and a holomorphic surjective map  $f:S\to B$  whose general fibers are all smooth elliptic (such a map is usually called {\em an elliptic pencil}). 

In fact, the situation is even more precise:  any non-K\"ahler compact surface with algebraic codimension one can be obtained starting from an elliptic principal bundle over a smooth projective curve by taking unramified quotients and blow-ups (see e. g. \cite[Lemma 1, Lemma 2]{_Brinzanescu:manuscripta_}; also  \cite[Proposition 3.17, Lemma 3.18]{_Brinzanescu:bundles_}). 

\hfill

\remark\label{_isotrivial_} In particular, all the smooth fibers of an elliptic pencil $f:S\to B$  on a non-K\"ahler compact surface $S$ with algebraic codimension one are isomorphic to some fixed elliptic curve $F$. Moreover, if all the fibers of $f$ are reduced, then all elliptic curves on $S$ are isomorphic to the general fiber $F$.

\hfill

It is  a natural question to ask if the above result could be generalized to higher dimensions. An immediate example to see that the verbatim extension is hopeless is by taking the product $X=M\times S$, where $M$ is a Moishezon non-K\"ahler manifold, and $S$ is an elliptic $K3$ surface of algebraic codimension one. The manifold $X$ is obviously non-K\"ahler and its algebraic codimension is one; yet there are infinitely many non-isomorphic elliptic curves on $X$.

However, there exists a class of non-K\"ahler manifolds for which a similar result can be proven in arbitrary dimensions, namely the class of {\em Locally Conformally K\"ahler} manifolds (LCK manifolds, for short; see  \ref{_def_LCK_}).  This class is wide enough:  for instance, it contains {\em all} non-K\"ahler surfaces with algebraic dimension one and many manifolds obtained starting with an arbitrary  algebraic cone of some projective space, by deleting the vertex and acting by a cyclic infinite group (see Subsection \ref{_LCK_basics_}).

\hfill

To state the  first result of the paper (\ref{_isotrivial_}), we need to recall the following notion.

\hfill

\definition\label{_isotrivial_def_} An {\bf isotrivial} family of curves is a surjective holomoprhic map $f:X\to B$ between complex manifolds, with $\mathrm{dim}(X)=\mathrm{dim}(B)+1$ such that all its smooth fibers are isomorphic.  

\hfill

\theorem\label{_isotrivial_} Let $M$ be a non-K\"ahler compact complex LCK manifold with algebraic codimension one. Then $M$ is bimeromorphic to a manifold $M^a$ which has an isotrivial family of elliptic curves $f^a:M^a\to B^a$ over a smooth projective manifold $B^a.$

{\bf Proof.} See Section \ref{_section_3_} \endproof

\hfill

Unfortunately,  it seems difficult to produce an elliptic bundle bimeromorphic to $M$, as in the case  of  surfaces. The main difficulty in the case when $\mathrm{dim} (M)\geq 3$ comes from the fact that $f^a$ may have fibers of dimensions at least two.

A first attempt in this direction was done in \cite{_Angella_Parton_Vuletescu_} where it was proven that any compact complex LCK threefold of algebraic codimension one  can be obtained  from elliptic principal bundles over projective smooth bases, up to alterations.

However, the main result of \cite{_Angella_Parton_Vuletescu_}, namely Theorem 2.12, is valid only under two assumptions, seemingly true, but still not proven up to now: (1) the Strong Factorization Conjecture, and (2) the existence of minimal smooth models for elliptic non-K\"ahler compact threefolds.  

In the present  paper, we succeed to avoid the use of the above two conjectural assumptions and clarify the structure of compact LCK manifolds of arbitrary dimension and algebraic codimension one.

\hfill 

To fix the terminology, we shall use the following terms:
\begin{itemize}
\item A {\em modification} (see e. g. \cite[Definition 2.1]{_Ueno_book_}) is  a proper, surjective, holomorphic bimeromorphic  map $f:X\to Y$ between complex analytic spaces.
\item An   {\em alteration} (see  \cite{_de_jong_}, \cite{_abramovich_oort_}) is  a proper, surjective, generically finite map $f:X\to Y$ between complex analytic spaces.  Notice that any modification is a particular case of alteration.
\end{itemize}

We can now give the the precise statement of the main result of the present paper.

\hfill

\maintheorem\label{_structure_thm_}
Let $M$  be a non-K\"ahler  LCK smooth compact manifold with algebraic codimension one.  Then there exists an alteration $X\to M$  (with $X$ smooth)  and a projective  manifold $B$ such that $X$ is a  modification of an elliptic principal bundle over $B$.

\hfill

The proof (see Section \ref{_Main_theorem_proof_}) is mainly based on algebraic methods in analytic spaces, and goes on the following lines. Through a proper modification of our LCK manifold $M$, one can  produce an algebraic reduction, $f^a\ :\ M^a\to B^a$, whose general fibers are smooth elliptic curves (\cite[Theorem 12.4]{_Ueno_book_}). Since not all the fibers of the algebraic reduction are of dimension one, we are forced to apply the Hironaka flattening (\cite{_Hironaka:flattening_})  to $f^a$. We then obtain a map  $f\ :\ X\to B$, where $X$  a complex space (possibly singular), which is a modification of $M^a$, $B$ is a smooth projective proper modification of $B$, and all fibers of $f$ are of dimension one. After passing, if necessarily, to finite covers of $X$ and $B$, we can assume that $f$ has reduced fibers. Blowing-up convenient subspaces of $X$, we can assume $X$ smooth. At this point, the LCK structure of $M$ becomes essential and leads to a clear description of the fibers of $f$: they all consists of a smooth elliptic curve plus some trees of rational curves. Eventually, we prove that the manifold $X$ is a proper modification of an elliptic principal bundle.

\section{Preliminaries}
In this section we shall briefly review several basic notions concerning the geometric and analytic aspects of the manifolds we are interested in.

\subsection{Locally Conformally K\"ahler geometry}\label{_LCK_basics_}

In this subsection we introduce the
basic notions of locally conformally K\"ahler
(LCK) geometry. For details and proofs we
refer to \cite{_OV:book_}.

\hfill

\definition\label{_def_LCK_}
A  Hermitian manifold 
$(M,I,g,\omega)$ is called 
{\bf  locally conformally K\"ahler} (LCK) if 
there exists a closed 1-form $\theta$ such that $d\omega=\theta\wedge\omega$.
The 1-form $\theta$ is called the {\bf Lee form}. 

\hfill

\remark  If the Lee form is exact, $\theta=df$, then the manifold is  {\bf Globally Conformally K\"ahler} (GCK), since the metric $e^{-f}g$ is K\"ahler.

\hfill

\proposition \label{_Covering_definition_}
A  Hermitian manifold 
$(M,I,g,\omega)$ is LCK if and only if 
it admits a K\"ahler cover $(\tilde
M,\tilde\omega)\arrow M$ whose deck group $\Gamma$
acts on the K\"ahler metric $\tilde\omega$ 
by holomorphic homotheties. 

\hfill

\remark Since all smooth compact curves are K\"ahler, the LCK geometry is not interesting in dimension one. The examples in dimensions greater than one are abundant (see \cite{_OV:book_}), e. g. : almost all compact complex surfaces, all Hopf manifolds (either linear or non-linear), Oeljeklaus-Toma manifolds, Kato manifolds.

\hfill

\theorem (Vaisman, \cite{_Vaisman_Trans_})\label{_Vaisman_theorem_} 
A compact LCK manifold with non-exact Lee form does not
admit any K\"ahler metric.

\hfill

We proved in \cite[Corollary 2.14]{ovv1} that the LCK class is not well-behaved with respect to blow-ups. Indeed, the blow up at points preserves the LCK class, whilst the blow-up along a submanifold is LCK if and only if the restriction of the Lee form to the submanifold is exact. This motivated a natural extension of the LCK definition, such as the larger class be closed to blow-up.

\hfill

\definition {(\cite{_Angella_Parton_Vuletescu_})} 
A {\bf weak locally conformally K\"ahler} (WLCK) structure
on a complex manifold $(M,I)$ is given by a (1, 1)-form $\omega$ and a real 1-form
$\theta$ such that:
\begin{itemize}
	\item $d\omega=\theta\wedge\omega$ and $d\theta=0$;
	\item $\omega$ is strictly positive definite outside a proper analytic subset (called the {\bf degeneracy locus} of $\omega$).
\end{itemize}

\hfill


\remark\label{_WLCK_on_gen_finite_Remark_} If  $f\ :\ X\to Y$ is an alteration between manifolds $X$ and $Y$, and if $Y$ is  WLCK, then $X$ is  WLCK too (and the Lee form on $X$ is the pull-back of the Lee form on $Y$). Notice also that if $X$ is a WLCK manifold and $Y\subset X$ is a submanifold such that $Y$ is not included in the degeneracy locus of the WLCK structure of $X$, then $Y$ inherits a WLCK structure too. Moreover, Vaisman's theorem \ref{_Vaisman_theorem_} holds in WLCK geometry too.

\hfill

The following is an useful tool, firstly proved in LCK geometry, but which subsists in WLCK geometry too.

\hfill

\proposition (\cite{opv,_Angella_Parton_Vuletescu_})\label{_Lemma_on_fibrations_} 
 (“Lemma on fibrations”). Let $X$ and $B$ be complex manifolds
with ${\rm dim} X > {\rm dim} B$. Let $f : X\to B$ be a surjective proper holomorphic map
with connected fibres. Let $(\omega, \theta)$ be a WLCK structure on $X$. If the Lee class
$[\theta]=f^*[\alpha]$ is  in the image of the pull-back $f^*:H^1(B)\to H^1(X)$ induced by $f$ in cohomology, 
then $[\theta] = 0$.

\subsection{Algebraic reduction of complex manifolds}

We begin by recalling the concept of {\em algebraic reduction} of a compact complex manifold $M$ (see \cite{_Ueno_book_} for more details). 
 
By a classical result of Siegel, the space  ${\mathcal M}(M)$ of global meromorphic functions on $M$ is a finitely generated $\C-$algebra.

\hfill

\definition The transcendence degree of  ${\mathcal M}(M)$ is called the {\bf algebraic dimension},  $a(M)$, of $M$ and the difference $\dim(M)-a(M)$ is called {\bf algebraic codimension}. 

\hfill

\condefinition\label{algred}  Let $V$  be the affine algebraic variety associated to ${\mathcal M}(M)$ (that is, its field of rational global functions is ${\mathcal M}(M)$) and   $\Phi\ :\ M\dashrightarrow V$ a meromorphic dominant map. Replacing $V$ by its projective closure $\overline{V}$ and possibly resolving the singularities of $\overline{V}$, we obtain a smooth projective manifold $B^a$ and a meromorphic dominant map $\varphi\ :\  M\dashrightarrow B^a.$ Resolving the indeterminacies of $\varphi$, produces a proper modification $c:M^a\to M$ of $M$ and a holomorphic map $f^a\ :\  M^a \to B^a.$ The triple $(M^a, f^a, B^a)$ is called  {\bf an algebraic reduction} of $M$ (some authors call $B^a$ alone the algebraic reduction). Notice  that the algebraic reduction is only defined up to bimeromorphisms. 

\hfill

\remark If $M$ has an LCK structure, then $M^a$  is a modification of $M$, and it has an WLCK structure. The   Lee form of this structure is the pull-back by $c$ of the Lee form on $M.$


\hfill

The algebraic reduction is the main tool used by Kodaira in the classification of non-projective compact complex surfaces. For non-K\"ahler surfaces $S$ with $a(S)=1$, the algebraic reduction implies that $S$ is bimeromorphic to an elliptic fibration $\hat{S}\arrow B$ onto some smooth projective curve $B$ whose smooth fibers are  {\em isomorphic} elliptic curves. Let $S_{\min}$ be the minimal model of $S$, then  any singular fiber of the map $S_{\min}\arrow B$ admits a  smooth reduction, which is still an elliptic curve (see e. g. \cite[Lemma 1]{_Brinzanescu:manuscripta_}). Taking a suitably chosen (ramified)  cover of $B'\arrow B,$  and letting $S'\to S$ the induced (\'etale) cover, one can prove that the induced map $f'\ :\  S'\arrow B'$ has only smooth fibers (see \cite[Lemma 2]{_Brinzanescu:manuscripta_}), which are  smooth isomorphic elliptic curves. It follows that $S'\arrow B'$ is a locally trivial elliptic bundle.  Eventually, the non-K\"ahlerianity of $S'$  combined with a Kodaira's projectivity criterion for surfaces  implies that $S'\arrow B'$ is an elliptic principal bundle (see e. g.  \cite{_ovv:surf_}). Summing up, one has the following  description of the analytic structure of any non-K\"ahler surface:

\hfill

\theorem\label{_surf_}
Let $S$ be a non-K\"ahler smooth compact surface with algebraic codimension one. Then there exists an alteration (in fact, even a unramfied cover) $\eta:S'\arrow S$ with $S'$ smooth,  and a projective smooth curve $B$ such that $S'$ is a blow-up of an elliptic principal bundle $X'\to B.$

\hfill

\remark Since every non-K\"ahler elliptic surface admits an LCK metric (see \cite{_Belgun_,_Tricceri_, _ovv:surf_}), in \ref{_surf_} we can replace  ``non-K\"ahler" by ``LCK".

\hfill

\subsection{Flat morphisms and flattening}
We recollect some basic facts about flat morphisms that will be used in the sequel. In what follows, $X,Y...$ will denote complex analytic spaces.

\hfill

\definition
Let $X,Y$ be analytic spaces, and $x\in X$. A map  $f\ :\ X\to Y$  is called {\bf flat at $x$} if the local ring ${\mathcal O}_{x, X}$ is a flat module over ${\mathcal O}_{f(x), Y}.$ The map $f$ is called {\bf flat} if it is flat at all points in $X.$ The set of points $y\in Y$ such that $f$ is flat in any point of the preimage $f^{-1}(y)$ is called {\bf the flat locus} of $f$.

\hfill

\remark It is not hard to show that flatness commutes vith base change, that is, for any flat map $f:X\to B$ and any  map $f':B'\to B$, the map  induced by base-change $f':X'=X\times_B B'\to B'$ is also flat.

\hfill

\remark Flat maps are {\bf equidimensional} (see e. g. \cite[Proposition 2.11]{_Grauert_Peternell_}, that is, for any point $y\in Y$ the fiber $X_y:=X\times_Y \{y\}$ has the same dimension.  {The converse does not hold in general. Still, an important result in this direction is \ref{_miracle_flatness_}. Before stating it, we need to  recall some basic notions of commutative  algebra.

\hfill

\definition\label{_CM_} (\cite{_bruns_herzog_}) Let $(A, \mathfrak{m})$ be a  local noetherian ring. A {\bf system of parameters} for $A$ is a sequence $t_1, \dots, t_n\in \mathfrak{m}$ such that $t_1$ is not a zero divisor in $A$ and for $i\geq 2$, $t_i$   is not a zero-divisor in $A/(t_1,\dots, t_{i-1}).$ The {\bf depth} of $A$ is the maximal length of a system of parameters. The ring $A$ is called {\bf Cohen-Macaulay} (CM,  for short) if its depth equals its Krull dimension, $\mathrm{depth}(A)=\dim(A).$} An analytic space $X$ is called CM if all its local rings ${\mathcal O}_{x, X}$ are CM.

\hfill

For a proof, see e. g. \cite[Theorem 2.2.4]{_bruns_herzog_}, (resp.  \cite[Theorem 2.1.3]{_bruns_herzog_}).


\hfill

\hfill

\lemma\label{_base_alg_} Let $A$ be a local noetherian ring and $t\in A$ not a zero divisor. If $A/(t)$ is regular (resp. Cohen-Macaulay) then $A$ is regular (resp. Cohen-Macaulay).
In particular, if $X_C$ is smooth (resp. Cohen-Macaulay) then $X$ is smooth (resp. Cohen-Macaulay) in all points of $X_C$.


\hfill

\proposition\label{_miracle_flatness_} {(\bf Miracle flatness}, \cite[Theorem 23.1]{_matsumura_})\label{MF} Let  $X,Y$ be analytic spaces, with $X$  Cohen-Macaulay and  $Y$ smooth. Then any equidimensional surjective morphism $f\ :\ X\to Y$  is flat.

\hfill


\remark A situation when the flatness of a morphism $f:X\to Y$ can be easily proven is when $X$ is locally irreducible of dimension $2$, $Y$ is smooth of dimension one and $f$  is a surjective morphism  $f:X\to Y$. In this case  for any $x\in X$ the ring $A={\mathcal O}_{x, X}$ is a domain and the ring $R={\mathcal O}_{f(x), Y}$ is a discrete valuation ring (DVR), since $Y$ is smooth. But this immediately implies that $f$ is flat, since as $A$ is a domain, it is a torsion-free $R-$module, and since $R$ is a DVR, torsion free $R-$modules  are flat (see e. g. \cite[Corollary 6.3]{_eisenbud_}. 

\hfill

The typical example of a non-flat map is the blow-up map, as it is non equidimensional. In particular, non-flat maps arise usually as soon as one does a resolution of singularities. It is thus natural to ask whether one can turn a non-flat map into a flat one, by performing modifications. The cornerstone result is due to Hironaka (first proven by M. Raynaud and L. Gruson, \cite{_Raynaud_Gruson_}, in the context of algebraic geometry):

\hfill

\theorem\label{Hir Thm} (\cite{_Hironaka:flattening_}, Hironaka flattening) Let $f:X\to Y$ be a proper surjective morphism of complex spaces. Then there exists modifications  $m_X:X'\to X$, $m_Y:Y'\to Y$ and a proper surjective holomorphic map $f': X'\to Y'$such that the following diagram is commutative
$$
\begin{CD}
X@<{m_X}<<X'\\
@VfVV @VVf'V\\
Y@<<{m_B}<Y'
\end{CD}
$$  and $f'$ is flat. 

\hfill

A more concrete description of the flattening is given  in \cite[Corollary 1.1]{_Hironaka:flattening_} which we quote verbatim:

\hfill

\theorem For $f:X\to Y$ of the theorem, there exists a projective
bimeromorphic morphism $\pi: Y'\to Y$ and a closed complex subspace $D$ of $Y$ such
that:

(i) $D$ is nowhere dense in $Y$, $\pi^{-1} (D)$ is nowhere dense in $Y'$, and $\pi$
induces an isomorphism $Y'\setminus \pi^{-1}(D)\simeq Y\setminus D$.

(ii) if $X'$ is the closure of $ X\times_Y Y'\setminus (\pi\circ f)^{-1}(D)$  in  $ X\times_Y Y'$ then $f'$
induces a flat morphism $X'\to Y'$.
In particular, the fibres $(f')^{-1}(y')$ with $y'\in Y'$  have the same dimension over
each connected component of $Y'$.

\hfill

\remark In  \ref{Hir Thm}, the space $X'$ is obtained from $X$ by performing a finite sequence of well-chosen  blow-ups. Unfortunately,  there is no control on the singularities of $X'$. Still, by the commutativity of the flattening with the base-change, and using the Hironaka resolution of singularities \cite{_Hironaka_resolution_}, one may assume that $Y'$ is smooth.

\hfill

\example (i) Let $Y$ be a manifold of dimension at least $2$, and $y\in Y$. Let $X={\rm Bl}_y(Y)$ and $\sigma:X\to Y$ the blow-up map. A flattening of $\sigma$ is obtained by taking $X'=Y'=X$ and $f'=id_{X}.$ 

(ii) A more illuminating example appears in \cite[Example 1]{_Hironaka:flattening_}. It is obtained by taking $Y$  to be a smooth surface, $C$ a smooth curve and $M=Y\times C$ with {$p:\ M\to Y$ being the projection on $Y$}.  Fix $y\in Y$ and let   $c\in \{y\}\times C\subset M$. Consider the blow-up $X={\rm Bl}_{c}(M)$  and let $f$ be the composition of $p$ with the blow-up map $X\to M.$ A flattening of $f$ is obtained by taking  $X'$ to be the blow-up of $X$ along the strict transform of the fiber $\{y\}\times C\subset M$ of $p$ in $X$, $B'={\rm Bl}_y(Y)$ and defining $f':X'\to Y'$ to be  the induced map.


\hfill

For further use we recall the following well-known-fact.

\hfill

\lemma\label{_base_alg_} Let $A$ be a local noetherian ring and $t\in A$ not a zero divisor. If $A/(t)$ is regular (resp. Cohen-Macaulay) then $A$ is regular (resp. Cohen-Macaulay).

\section{Proof of \ref{_isotrivial_}}\label{_section_3_}

Let $f^a:M^a\to B^a$ be an algebraic reduction and $c:M^a\to M$ a modification as in \ref{algred}. Let $E_c\subset X$ be the exceptional divisor of this modification. Let $Z\subset B^a$ be the {\bf non-flat locus} of $f$, that is, $Z=\{b\in B^a\mid \mathrm {dim}((f^a)^{-1}(b))\geq 2\}$   it is  a closed analytic subspace of $B^a$   (\cite{_frisch_}, \cite[3.18]{_fischer_}) of codimension at least 2. By \cite{_Ueno_book_}, the general fibers of $f^a$ are smooth elliptic. By our assumption on $M$, there exists an LCK metric $\omega$ on it with associated Lee form $\theta.$ This induces a WLCK structure on $M^a$ with Lee form $\theta_{M^a}:=c^*(\theta)$; denote by $D_\omega\subset M^a$ the degeneracy locus of this WLCK structure.   Choose a smooth curve $C\subset B^a$. Since $Z$ is of codimension at least $2$,  we may assume that $C\cap Z=\emptyset;$ then $M^a_C:=M^a\times_B C$ is a surface. Choosing $C$ general, we may assume that the general fibers of $f_C: M^a_C\to C$ are smooth (where $f_C=f^a\restrict{M^a_C}$) and that $M^a_C$ is not contained in $D\cup E_c$.  Then $c(M^a_C)\subset M$ is a surface, because  $M^a_C\not\subset E_c$. Let now  $\sigma: \widehat{M^a_C}\to M^a_C$ be a desingularization of $M^a_C$, and observe that  since $M^a_C\not \subset D_\omega$,  $\widehat{M^a_C}$  has  a WLCK structure, induced by the pull-back via $c\circ \sigma$, with Lee form $\theta_C:= (c\circ \sigma)^*(\theta)$.   
	
We now prove that the WLCK structure on $\widehat{M^a_C}$ cannot be K\"ahler. Assume, by contradiction, that $\widehat{M^a_C}$ is K\"ahler. Then, Vaisman's theorem for WLCK structures (\cite{_Angella_Parton_Vuletescu_}) implies that $[\theta_C]=0.$  Since $\sigma$ induces an injection at the $H^1$-level, we deduce that $[\theta_{M^a}]\restrict{ M^a_C}=0.$ Therefore, as $C$ is general,  the restrictrion of $[\theta_{M^a}]$ to the general  fibers of $f$ vanishes,  and hence $[\theta_{M^a}]$ is the  pullback of a cohomology class form $B^a$. Now \ref{_Lemma_on_fibrations_}  implies that $[\theta_{M^a}]=0.$  Finally, we note that the natural  map $H^1(M)\to H^1(M^a)$ is an isomorphism, because $M$ and $M^a$ are bimeromorphic< therefore, $[\theta]=0$, that is $M$ is  K\"ahler, absurd. 

Eventually, we prove that all smooth fibers of $f$ are isomorphic. Indeed, the general fibers of $f$ are the general fibers  of the map $M^a_C\to C$, hence also of the map  $\widehat{M^a_C}\to C$.  Since $\widehat{M^a_C}$ is non-K\"ahler for general $C$, it these fibers are isomorphic.  Then the continuity of the $j-$invariant implies that all smooth fibers of $f$ are isomorphic, hence $f:M^a\to B$ is an isotrivial family as stated. \endproof

\section{Proof of the  \ref{_structure_thm_}}\label{_Main_theorem_proof_}

\subsection{A technical result}

The following is the main technical result we need.

\hfill

\theorem\label{_smooth_flatt_} Let $M$ be an LCK manifold, and $f^a:M^a\to B^a$ an algebraic reduction of $M$ with $B^a$ smooth projective. Then there exist the manifolds $X$ and $B$, with $B$ smooth projective, a flat  morphism $f:X\to B$, and the alterations $r_X:X\to M^a$, $r_B:B\to B^a$ such that the diagram
$$
\begin{CD}
M^a@<{r_X}<<X\\
@Vf^aVV @VVfV\\
B^a@<<{r_B}<B
\end{CD}
$$
is commutative.

\hfill

\proof 
We describe  first the idea of our argument.

By \cite{_Hironaka:flattening_}, there exist a compact complex space $X$, a projective manifold $B$, a flat map $f:X\to B$  and the modifications 
 $m_X:X\to M^a$, $m_B:B\to B^a$ such that the diagram
$$
\begin{CD}
M^a@<{m_X}<<X\\
@Vf^aVV @VVfV\\
B^a@<<{m_B}<B
\end{CD}
$$
is commutative.

The issue is that Hironaka flattening does not guarantee that $X$ is smooth. To get a smooth $X$, we procced as follows.

In a first step, in \ref{_reduction_} we perform alterations  of $X$ and $B$ to eliminate the  possible non-reduced fibers. As a result, we obtain a diagram
\begin{equation}\label {_reduc_}
\begin{CD}
X@<{r_X}<<X'\\
@VfVV @VVf'V\\
B@<<{r_B}<B'
\end{CD}
\end{equation}
with $B'$ smooth projective and $f'$ flat with all its fibres reduced.

In a second step, we use the LCK structure of $M$ to show (\ref{_nk_})  that  any fiber of $f'$ consists of an elliptic curve plus some trees of rational curves (this will be also useful in the last part of the paper).

Eventually, using \ref{_blow_up_}, we show that there exist the alterations  $r_{X'}:\hat{X}\to X'$ and  $r_{B'}:\hat{B}\to B'$ with both $\hat{X}$ and $\hat{B}$ smooth, and an equidimensional map $\hat{f}: \hat{X}\to \hat{B}$ such that the diagram
$$
\begin{CD}
X'@<{r_X'}<<\hat{X}\\
@Vf'VV @VV\hat{f}V\\
B'@<<{r_B'}<\hat{B}
\end{CD}
$$
is commutative.

\hfill

We now provide the details of the proof.

We start with the following lemma wich is a weak version of the ``Weak semi-stable reduction" theorem (see e. g. \cite{_Abramovitch_Karu_}, \cite{_de_jong_}). Here we only consider the case of equidimensional families of curves, and we  only ask that the resulting fibers be reduced (we do not pay attention to their semi-stability). For convenience of the reader, we include a proof for our context.

\hfill

\lemma\label{_reduction_} Let $f:X\to B$ be  a flat surjective  morphism, with $B$ smooth projective and  ${\rm dim}(X)={\rm dim}(B)+1$. Then there exist  alterations $r_X:X'\to X$, and  $r_B:B'\to B$ (with $B$ smooth projective) and a surjective flat holomorphic map $f':X'\to B'$ such that the diagram

$$
\begin{CD}
X@<{r_X}<<X'\\
@VfVV @VVf'V\\
B@<<{r_B}<B'
\end{CD}
$$
is commutative and such that all the fibers of $f'$ are reduced.

\hfill

\proof  Let $Z\subset X$ be the locus of the fibers with nilpotent elements and $D=f(Z).$ Using Hironaka resolution of singularities,  there exists a composition of blow-ups $c: \hat{B}\to B$ with $\hat{B}$ smooth and such that the inverse image $\hat{D}\subset \hat{B}$ of $D$ is a divisor in simple normal crossing. Let $\hat{D}=\sum_i m_iC_i$ be the decomposition of $\hat{D}$ into irreducible components. We apply the Kawamata trick (\cite[Theorem 17]{_Kawamata_}) to $\hat{B}$ with respect to  $(m_i, C_i)$. This way, we produce a new manifold $B'$ and a covering $r:B'\to \hat{B}$, ramified of order $m_i$ along each $C_i$. Let $\hat X:=X\times_B \hat B$, and  $X'$  the normalization of $\hat{X}\times_{\hat{B}} B'$. With a computation similar to the one in \cite[Lemma 1.4]{_Angella_Parton_Vuletescu_}, one can see that the induced map $f':X'\to B'$ has only reduced fibers. Notice that the fibers of $f'$ are curves, since the normalization map $X'\to \hat{X}\times_{\hat{B}} B'$ is finite and the fibers of  $\hat{X}\times_{\hat{B}} B'\to B'$ are curves, and hence $f'$ is flat.  \endproof


\hfill

To proceed, we need to introduce the followibg definition.

\hfill

\definition\label{_slice_} Let $f:X\to B$ be a morphism of complex spaces and let $C\subset B$ be a curve. The {\bf slice} of $f$  along $C$ is the subspace $X_C=X\times_B C$ of $X$, that is $\mid X_C\mid  =f^{-1}(C)$ with ideal  sheaf ${\mathcal J}_{X_C\mid X}=f^*\left({\mathcal J}_{C\mid B}\right)$, where ${\mathcal J}_{C\mid B}$ is the sheaf of ideals of $C $ in $B.$

\hfill

\lemma\label{_nk_} Let $f':X'\to B'$ be a flat surjective morphism  with reduced fibers with $B'$ a projective manifold and ${\rm dim}(X')={\rm dim}(B')+1.$ Assume that the general fibers of $f'$ are smooth and connected. Then,  for a general curve $C\subset B$ we have:

(i) All the fibers of $f'$ are connected.

(ii) The slice $X'_C$ is normal. 

(iii)  If $X'$ is an alteration of an LCK manifold, then any desingularization $\widehat{X'_C}$ of  the slice ${X'_C}$ is a non-K\"ahler (elliptic) surface. 

(iv) If $X'$ is an alteration of an LCK manifold, any fiber of $f'$ consists of a smooth elliptic curve plus some trees of rational curves. Moreover, all elliptic curves contained in the fibers of $f'$ are isomorphic (to a fixed elliptic curve $F$).

\hfill

\proof 
(i) Consider the Stein factorization of $f'=h\circ f$ given by $X'\stackrel{f}{\longrightarrow} Y\stackrel{h}{\longrightarrow}B'$ with $f$ a connected morphism and $h$ a finite morphism. Since the general fibers of $f'$ are connected,  $h^{-1}(b)$ is a single point for general $b\in B'.$ Therefore $h$ is a bimeromorphism, and hence by Zariski Main Theorem $h^{-1}(b)$ is connected for all $b\in B'.$ The finiteness of $h$ implies  that $h^{-1}(b)$ is  a single point  for all $b\in B,$ and thus $h$ is an isomorphism (because $B'$  is normal).

(ii) First, we notice that $X'_C$ is irreducible: indeed, $C$ is smooth and connected (being general), and also the general fiber of the fibration $X'_C\longrightarrow C$ is smooth and connected. 
Next, we notice that since the fibers of the induced map $f\restrict{X'_C}:X'_C\to C$ are  one-dimensional and reduced, they are Cohen-Macaulay (see e. g. \cite{_Burban_Drozd_}). Therefore, \ref{_base_alg_} implies that $X'_C$ is CM.  As  there are finitely many singular points in the fibers, using again \ref{_base_alg_} we deduce that $X'_C$ has singularities at most at the points where the fibers are singular. In particular, $X'_C$ has finitely many singular points. The normality of $X'_C$ now follows from  
Serre's criterion for normality  (see e. g \cite[Theorem 23.8]{_matsumura_}) which in our context reads:   ``An analytic space of dimension $2$ is normal if and only if it is Cohen-Macaulay and it has at most finitely many singular points". 
This proves  (ii).

(iii) We denote by $\sigma:\hat{X}'_C\longrightarrow X'_C$  a desingularization of $X'_C.$
Let $S_C\subset M$ be the image of $X'_C$ under the composition of alterations $X'\longrightarrow X\longrightarrow M^a\longrightarrow M$. Since $C$ is general,  $S_C$ is a surface. We thus obtain a proper, surjective, generically finite morphism  $s:\hat{X}'_C\longrightarrow S_C$. Pulling-back the LCK structure from $M$ to $\hat{X}'_C$ we see that $\hat{X}'_C$ has a WLCK structure with Lee form $s^*(\theta)$ (\ref{_WLCK_on_gen_finite_Remark_}).

Assume, by absurd, that $\hat{X'}_C$ is K\"ahler. Recall that Vaisman theorem also holds for WLCK structures (\ref{_WLCK_on_gen_finite_Remark_}), therefore the class $s^*([\theta])=0.$ This implies that the restriction of $[\theta_{X'}]$  to the general fibers of  $X'_C\longrightarrow C$ is also zero. Since $C$ is general, by moving it we can cover the whole $B$, and hence  the restriction of  $[\theta_{X'}]$ to the general fibers of  $X\longrightarrow B$ vanishes. It follows that the restriction of $[\theta_{X'}]$ to the preimage of the set $\Reg(f)$ of regular values of $f$  is the  pull-back of a class  $\alpha\in H^1(\Reg(f))$. This class $\alpha$ can be seen as a cohomology class on $B$, since  $H^1(B)$ and $H^1(\Reg(f))$ are canonically isomorphic.

We conclude that the restriction of $\alpha$ to  iterated hyperplane sections vanishes, and hence, by  Lefschetz, $\alpha=0$. Therefore  $[\theta_{X'}]=0$ which implies $[\theta]=0,$ absurd, hence $\hat{X}'_C$ is non-K\"ahler.  This proves  (iii).

\hfill

 The statement (iv) follows from (ii) and  (iii). Indeed, all the fibers of $f'$ are also fibers of general slices $X'_C\to C$. On the other hand, $X'_C$ is normal, and hence the fibers of $X'_C\to   C$ are obtained by contracting some curves in the fibers of the desingularization $f'\restrict{X'_C}\circ \sigma : \hat{X}'_C\to C.$ Since $\hat{X}'_C$ is non-K\"ahler, its fibers are smooth elliptic curve plus trees of rational curves \cite{_ovv:surf_}. 

The elliptic component is not contracted in any fiber $\left(X'_C\right)_c=({f'}\restrict{ {X'_C}})^{-1}(c)$ ($c\in C$). Indeed, since contracting rational curves from a tree yields a tree,  such a  singular fiber  would be simply-connected. This implies that there exists some open neiborhood $V$ of $\left(X'_C\right)_c$ which is simply-connected. It follows the the restriction of the class $[\theta_{X'}]$ of the Lee form to any  fiber of $X'_C\to C$ contained in $V$  vanishes, hence the restriction of  $[\theta_{X'}]$ to the general fibers of $X_C\to C$ vanishes, which implies that $[\theta_{X'}]$ is a pull-back from $C$. But then $\hat{X}'_C$ would be K\"ahler by  \ref{_Lemma_on_fibrations_}.
Eventually, the fact that all elliptic curves contained in the  fibers are isomorphic  follows from the fact that this is true for the fibers of any non-K\"ahler elliptic surfaces. This proves (iv). \endproof

\hfill

For the last step of the proof of \ref{_smooth_flatt_}, we introduce a numerical invariant as a measure of singularities of $X.$

\hfill

\definition  Let $f:X\to B$ be a flat surjective morphism with reduced fibers, with $B$ a projective manifold, and ${\rm dim}(X)={\rm dim}(B)+1$. Assume that the general fibers of $f$ are smooth.

 (i) Fix an arbitrary point $b\in B$. For an arbitrary general curve $C\subset B$ through $b$ we define
{\bf the defect} ${\Def}(X_C)$ of $X_C$ as the number of irreducible components of the fiber $F_b$ which are contracted (i. e. mapped to points) to $X_b$ under the minimal desingularization  $\sigma:\hat{X}_C\longrightarrow X_C$ of $\hat{X}_C$. 

(ii) {\bf The defect} $\Def(b)$ of a point $b\in B$  is the minimum of ${\Def}(X_C)$,  where $C\subset B$ is a general curve through $b.$ 

(iii) {\bf The defect} $\Def(X)$ of $X$ is defined as $\Def(X):=\max_{b\in B}\Def(b).$

\hfill

\remark\label{_X_smooth_along_X_b_cor_} Suppose ${\Def}(X_C)=0$.   Then $\sigma$ is an isomorphism, hence ${X}_C$ is smooth. It follows that $X$ is smooth along $X_C$ (see \ref{_base_alg_}). 
In particular, if $b\in B$ is a point such that if ${\rm Def}(b)=0$, then $X$ is smooth along $X_b.$ 

\hfill

\lemma\label{_blow_up_}  Let $f':X'\to B'$ be a flat surjective morphism with reduced fibers with $B'$ a projective manifold and ${\rm dim}(X')={\rm dim}(B')+1$.  Assume that the general fibers of $f'$ are smooth. Then there exists a modification $m_X\ :\ X''\to X'$, a modification $m_B\ :\ B''\to B'$ with $B''$ a smooth projective manifold, and a flat holomorphic map $f'':X''\longrightarrow B''$ such that the following diagram is  commutative
$$
\begin{CD}
	X' @<m_X<< X'' \\
	@V f' VV @VV f'' V \\
	B' @<<m_B< B''
\end{CD}
$$
and such that ${\Def}(X'')<{\Def}(X').$

\hfill

\proof 
Let $S(X')$ be the union of singular loci of the fibers of $f',$ 
$$S(X')=\bigcup_{b\in B} \mathrm{Sing}(X'_b).$$
Since any fiber has at most finitely singular points, we obtain that  ${\rm dim}S(X')={\rm dim}f'(S(X'))$.
Since  the general fibers of $f':X'\to B'$ are smooth, we see that $f'(S(X'))$ is a proper subspace of $B'$, hence ${\rm dim}f'(S(X'))<{\rm dim}(B')$.
If ${\rm codim}(f(S(X'))\geq 2$,  by blowing-up $f'(S(X'))$ and possibly desingularizing ${\rm Bl}_ {f'(S(X'))}B'$ (using \cite{_Hironaka_resolution_}) we get a modification $b_B:B_1\to B'$. Doing base change: 
$$
\begin{CD}
	X' @<b_X<< X_1 \\
	@V f' VV @VV f_1 V \\
	B' @<<b_B< B_1
\end{CD}
$$
where  $X_1=X'\times_{B}B_1$, we may assume that  $S(X_1)$ is of codimension $2$.
Notice that since we are doing  base-change, the morphism $f_1$ is still flat and,  for  any point $b\in B_1$,  the fibers of $f_1$  are isomorphic to the  fibers of  $f$. In particular, ${\rm Def}(X')={\rm Def}(X_1).$

Let $\sigma:\overline{X}\to X_1$ be the blow-up of $X_1$ along $S(X_1)$ (considered with the reduced structure). Since $S(X_1)$ is of codimension $2$,  then for any $b\in B_1$ the fiber of the map $\overline{f}=f_1\circ \sigma:\overline{X}\to B_1$ over  $b$ is isomorphic to  the fiber of the map $f_1$ over $b$ plus some newly created curve which is the fiber of the natural map from the exceptional divisor of $\sigma$ to   $S(X_1)$ in $x$. Since  the defect counts the number of contracted curves, this construction decreases the defect. Note that blowing-up along a reduced space does not create nilpotent elements, and hence all the fibers of $\bar f$ are  reduced. Therefore ${\rm Def}(\overline{X})<{\rm Def}(X_1)={\rm Def}(X).$  
\endproof

\hfill

This comppletes the  proof of \ref{_smooth_flatt_}.  \endproof


\subsection{The proof of the main theorem}

Let $M$ be LCK. Using the algebraic reduction (\ref{algred}), there exists a modification $M^a\to M$ such that $M^a$ sends a surjective holomorphic map $f^a: M^a\to B^a$ onto a projective smooth manifold.  Using Hironaka flattening (\ref{Hir Thm}), there exists a modification $X_1$ of $M^a$ and a flat morphism $f:X_1\to B_1$ onto some smooth projective manifold $B_1.$ Using \ref{_reduction_}, we produce an alteration $X_2\to X_1$ such that there exists a morphism $X_2\to B_2$ as in diagram \eqref{_reduc_}  having only reduced  fibers. Making repeated use of \ref{_blow_up_}, we produce a new modification $X$ of  $X_2$ which is smooth and  a flat map
$f:X\to B$ onto a projective manifold such that all its fibers $X_b$ are composed of an elliptic curve ${\mathcal E}_b$ plus some trees of rational curves.  Notice that there exists a smooth non-K\"ahler  slice $X_C$ through every point of $X$ because ${\rm Def}(X)=0$.

\hfill

Eventually, to produce an elliptic bundle ${\mathcal P}\to B$ such that $X$ is a modification of ${\mathcal P}$, we argue as follows.

Recall the context. We have compact complex manifolds $X$ and $B$, a fixed smooth curve $F$ of genus one, and a holomorphic surjective map $f:X\to B$  such that, any fiber $X_b$  of $f$ consists of a genus one curve ${\mathcal E}_b$   isomorphic to $F$ plus some trees $\{T^i_{b}\}_{i=1,\ldots m(b)}$  of rational  curves attached to ${\mathcal E}_b$ at the  points $R^i_b, {i=1,\dots m(b)}$. We  call ${\mathcal E}_b$ {\bf the stem of the fiber}, and  $R^i=R_{T^i}$  the {\bf root points} of the fiber $X_b$. Notice that each stem ${\mathcal E}_b$ is isomorphic to a fixed elliptic curve $F$.  Moreover,  for any point $b\in B$ such that the fiber $X_b$ is singular, there exists a smooth curve $C\subset B$ such that its preimage $X_C=X\times_BC$ is a smooth non-K\"ahler surface.
We shall also denote by $\Delta\subset B$ the {\bf discriminant} of $f$, that is,  the set  of critical values of $f$. By Sard Lemma, $\Delta$
is a proper subvariety of $B$. Let ${\mathcal R}\subset X$  be the set of all root points contained in the fibers of $f$; one can that ${\mathcal R}$ is of codimension  at least $2.$

Choose a local section $s: U\to X$ of $f:X\to B$. For an appropriate choice of $s$,  its image $\Sigma=s(U)$ contains no root point of any fiber. Let $D=f^{-1}(U)$ and consider  the linear map $\varphi:D\to {\mathbb P}(f_*\left({\mathcal O}_D(3\Sigma)\right)$; call $P$ the image of $\varphi$ (see the diagram below).
\[
\begin{tikzcd}[row sep=small, column sep=small]
	D \arrow[rr, "\varphi"] \arrow[rd, "f"'] & & P \arrow[ld, "p"] \\
	& U & 
\end{tikzcd}
\]
For any point $b\in U$, the restriction of the map $\varphi$ to the fiber $D_b=X_b$ contracts all the rational curves and maps isomorphically the stem ${\mathcal E}_b $ of the fiber to the fiber $P_b$ of $P\to U$ at $b$.
In particular, all the fibers of the map $P\to U$ are isomorphic to $F$.

We will show  that $P$ is smooth. Choose any point $b\in U;$ by the assumption on $f$, there exists a smooth curve $C$ in $U$ such that $X_C=X\times_U C$ is a smooth  surface. Let $P_C=P\times_U C$; then the restriction of  $\varphi:X_C\to P_C$ is the blow-down of all rational curves, hence $P_C$ is a smooth surface too. Using \ref{_base_alg_}  and induction it follows that $P$ is smooth. 

Hence, $p:P\to U$ is a fibration by smooth isomorphic curves: by the theorem of Fischer-Grauert (\cite{_Fischer_Grauert_}), it is a locally trivial $F$-bundle. Shrinking $U$ if necessary, we may assume it to be a principal bundle. 

Clearly, if one takes a different section $s'$ as above, considering the map $\varphi'$ and  $P'$, the  principal $F-$bundle associated as above, then $P\setminus\varphi({\mathcal R}\cap D)$ and $P'\setminus\varphi'({\mathcal R}\cap D)$ are isomorphic (via $\varphi\circ {{\varphi}}'^{-1}$). By Riemann second  extension theorem (\cite[p. 132]{Grauert-Remmert_coherent _analytic_sheaves})  we get an isomorphism of $F$-bundles $g:P\to P'$ (by passing to the universal covers)  .

\hfill




Now cover $B$ by open subsets $\{U_\alpha\}_{\alpha\in I}$ such that for each $U_\alpha$ we fixed a section of $f$, say  $s_\alpha: U_\alpha \to X$, not passing through any root points in  the fibers. Using the above procedure, on the overlappings $U_\alpha\cap U_\beta$  we obtain a section $g_{\alpha\beta}\in\Gamma\left(U_\alpha\cap U_\beta,   {\mathcal Aut}_X(F)\right)$. By their definitions, the $g_{\alpha\beta}'$s obey the cocycle relation, hence they define an element in the set $\breve H^1(X,  {\mathcal Aut}_X(F))$, that is, a locally trivial $F-$bundle $P$ over $B.$

Eventually, gluing  the maps  $\varphi_\alpha: f^{-1}(U_\alpha)\to P_\alpha$ together we get a map  $\varphi:X\to P$ which contracts all the rational curves on $X$.
\endproof


\hfill

{\bf Acknowledgment:} We are indebted to Marian Aprodu for countless discussions  which improved our proofs.

{\small

\noindent {\sc Liviu Ornea\\
	University of Bucharest, Faculty of Mathematics and Informatics, \\14
	Academiei str., 70109 Bucharest, Romania}, and:\\
{\sc Institute of Mathematics ``Simion Stoilow" of the Romanian
	Academy,\\
	21, Calea Grivitei Str.
	010702-Bucharest, Romania\\
	\tt lornea@fmi.unibuc.ro,   liviu.ornea@imar.ro}

\hfill

\noindent {\sc Misha Verbitsky\\
	{\sc Instituto Nacional de Matem\'atica Pura e
	Aplicada (IMPA) \\ Estrada Dona Castorina, 110\\
	Jardim Bot\^anico, CEP 22460-320\\
	Rio de Janeiro, RJ - Brasil }\\
	also:\\
	Laboratory of Algebraic Geometry, \\
	Faculty of Mathematics, National Research University 
	HSE,\\
	6 Usacheva Str. Moscow, Russia}\\
\tt verbit@verbit.ru, verbit@impa.br 

\hfill

\noindent {\sc Victor Vuletescu\\
	University of Bucharest, Faculty of Mathematics and Informatics, \\14
	Academiei str., 70109 Bucharest, Romania}\\
\tt vuli@fmi.unibuc.ro}

\end{document}